\documentclass[letterpaper, 10pt, conference]{ieeeconf}
\IEEEoverridecommandlockouts

\usepackage{xcolor}

\definecolor{subsectioncolor}{rgb}{0,0,0}
\usepackage{generic}
\usepackage{cite}
\usepackage{wrapfig}
\usepackage{amsmath,amssymb,amsfonts}
\usepackage{algorithmic}
\usepackage{graphicx}
\usepackage{algorithm,algorithmic}
\usepackage{hyperref}
\hypersetup{hidelinks=true}
\usepackage{textcomp}
\usepackage{dsfont}
\algsetup{indent=1.5em}
\usepackage{bm}
\usepackage{mathrsfs}
\usepackage[utf8]{inputenc}
\usepackage{nicefrac}

\usepackage{bm}
\usepackage{mathptmx}
\usepackage{tcolorbox}
\usepackage{tikz}
\let\labelindent\relax
\usepackage{amssymb, amsmath}
\usepackage{mdframed}
\usepackage[shortlabels]{enumitem}

\newcommand{\x}{{ {\bf x}}}

 \renewcommand{\S}{{\mathcal S}}

\newcommand{\p}{{\bf p}}
\newcommand{\bp}{{\bf p}}

\renewcommand{\u}{{ {\bf u}}}
\newcommand{\z}{{ {\bf z}}}
\newcommand{\eop}{{\hfill $\Box$}}
\newtheorem{thm}{{\bf Theorem}}

\newcommand{\ignore}[1]{}

\IEEEoverridecommandlockouts                              

\usepackage{bbm}
\usepackage{enumerate}
\newcommand{\TR}{\ifthenelse{1>2}}

\newcommand{\ForTR}[1]{}  
\usepackage{multirow}
\newcommand{\exactproof}[1]{}  

\newcommand{\Forarxiv}[1]{#1}   
\newcommand{\Forecc}[1]{} 

\newtheorem{propsn}{\bf {Proposition}} 

\newcommand{\muL}{\mu_{leb}}

 \newcommand{\Removesmoothapprox}[1]{}  
 \newcommand{\Removegeneralcost}[1]{}  

\newcommand{\balpha}{\bm{\alpha}}

\newcommand{\bbeta}{\bm{\beta}}

\begin{document}
\title{\bf Inventory Control and  Storage Design   via  an $L^{\infty}$-control approach  
}

\author{ Madhu Dhiman, Veeraruna Kavitha and Nandyala Hemachandra 
\thanks{Authors are with the Indian Institute of Technology Bombay. Emails: \{madhu.dhiman,vkavitha, nh\}@iitb.ac.in. }
}
\maketitle
\begin{abstract}
Inventory control traditionally focuses on controlling prices and production rates. However, such systems offer limited flexibility in adjusting physical storage capacity, which is typically fixed at the initial design stage. This paper addresses this important yet relatively unexplored aspect, leading to a non-standard optimal control problem. Specifically, the objective is formulated as a weighted combination of the classical integral term, capturing conventional inventory-related costs, and an $L^\infty$-term, representing the maximum inventory level over the entire planning horizon. We first establish a relationship between the optimal pricing and production-rate policies, thereby reducing the dimensionality of the problem.

Building on the framework developed in our earlier work, we establish the existence of an optimal relaxed control and, through a smooth approximation procedure, obtain an $\epsilon$-optimal pure control for the reduced problem. By numerically solving the resulting smooth approximate problem, we identify a policy that achieves a substantial reduction in the required storage capacity. Interestingly, this reduction entails only a negligible loss in net revenue (less than $6\%$), after accounting for shortage, holding, and production costs. Furthermore, we observe that increasing the weight assigned to the $L^\infty$-component leads to reduced fluctuations in the resulting inventory-level trajectory.

These findings also have important implications when information about future demand is incomplete or unavailable. In such settings, continually adjusting the pricing and production policies to keep inventory levels close to zero can substantially reduce storage requirements while incurring only a negligible loss in net revenue. Thus, our results highlight the potential of incorporating inventory-level considerations directly into operational decision-making, particularly when there is  flexibility to design the size of physical storage. \Forarxiv{We theoretically establish  certain interesting properties of the optimal  policy, which also support the above guidance.
}

\end{abstract}

\section{Introduction}
\label{sec_intro}

Traditional inventory control problems focus on optimizing shortage costs (arising from lost sales due to inventory unavailability), holding costs (associated with maintaining excess inventory), production costs (arising from time-varying production rates, which may differ from the rate best suited to the system), and other related costs, while catering to time-varying demands (see \cite{chen2004coordinating, chen2012pricing, elmaghraby2003dynamic, whitin1955inventory}).

An important yet relatively unexplored aspect is the physical storage capacity required to accommodate inventory. Unlike operational decisions such as pricing and production rates, physical storage capacity, once constructed, cannot be easily altered. Existing strands of literature consider constrained inventory optimization problems in which an upper bound is imposed on the maximum inventory level (see, e.g., \cite{tu2013introductory, zaher2013optimal}). These formulations are tailored to scenarios in which the storage infrastructure is already in place and its capacity is predetermined. However, to the best of our knowledge, the problem of optimally determining the required storage capacity, given a planning horizon and anticipated variations in future demand, has not been addressed in the literature. Such a design consideration is important in settings where there is flexibility in designing physical storage capacity, as it allows storage requirements to be optimized alongside conventional inventory-related costs.

It is important here to distinguish between inventory control within an existing infrastructure, where physical storage capacity is predetermined and imposed as a constraint (e.g., through barrier-function methods; see \cite{polak1993method, schiela2009barrier}), and the optimal design of the infrastructure itself, where storage requirements are treated as part of the optimization problem. These are fundamentally different problems and the  present paper addresses the latter.

Technically, the  problem discussed above leads to a non-standard control problem in which one aims to optimize (assume maximization for  the sake of convenience) a weighted combination of an  integral term  and an $L^{\infty}$ term, like   $\sup_u\{ \int_0^{t_1} L_t dt - \rho \sup_{t \in [0,t_1]} |F_t | \}$ (here $\rho$ is the trade-off factor);  for example, 
the revenue generated over the entire planning horizon considering holding and production costs leads to the integral term ($L_t$ at time $t$ represents the instantaneous combined cost), while the maximum inventory-level attained over the entire planning horizon leads to the $L^\infty$ term ($F_t$ is the inventory level at time $t$).

 There is a vast literature that studies classical control problems involving only integral terms like $\sup_u\{ \int_0^{t_1} L_t dt \}$ (e.g., \cite{li2015joint}), and a considerable size of literature that studies the  problems with only $L^\infty$ term, $\sup_u \{\sup_{t \in [0,t_1]} |F_t |\}$ (e.g., \cite{barron1989bellman}). To the best of our knowledge,  the literature combines the two costs in a restrictive  manner  and considers optimizing costs like,  $\sup_u \{ \sup_{t \in [0,t_1]} \{ \int_t^{t_1} L_s ds +   |F_t | \} \}$. The $L^\infty$ Bellman equation  supports only this type of combined objective functions (see \cite{barron1999viscosity}). Clearly, the combined  problem  modeling  the inventory control   is drastically different.

We address this non-standard optimal control problem by employing the framework developed in our recent  work \cite{dhiman2026optimal}. We introduce an auxiliary state variable, \(y(s):=\sup_{t\leq s}|F_t|\), to capture the instantaneous \(L^\infty\) term (e.g., it represents the maximum inventory level reached till time \(s\)) and convert the \(L^\infty\) term into a terminal cost. We establish the existence of an optimal relaxed control for the resulting problem. However, the dynamics of the auxiliary state variable involve an indicator function of the state variables, which makes the resulting control problem technically challenging. To address this difficulty, following the approach in \cite{dhiman2026optimal}, we consider an appropriate smooth approximation of the reformulated problem and establish the existence of an \(\epsilon\)-optimal pure control for the original problem.

The smooth approximate problem falls within the standard optimal control framework, allowing us to employ the Pontryagin Maximum Principle to derive the necessary optimality conditions and develop a computationally tractable numerical approach. Using this, we conduct an extensive numerical study of the inventory control problem and obtain several interesting observations.

As the importance of the \(L^\infty\) term increases, the variations in the optimal inventory-level trajectory diminish; however, the revenue at optimality does not reduce much (at most 6\% losses for highly sensitive markets and negligible losses for less sensitive markets). The theoretical and numerical observations culminate in a strategy advice to the manufacturer with zero future-demand information: \textit{it is beneficial to maintain the instantaneous inventory level near zero and to have a small shortage towards the end.} \Forecc{In our longer arxiv version \cite[Theorem 3]{dhiman2024optimal}, we  in fact prove a structural result to this effect, under an additional assumption. }
\Forarxiv{We  in fact prove a structural result to this effect, under an additional assumption. }

There are few strands of  inventory control literature that  consider a constraint optimization problem with a hard constraint on the inventory levels (e.g., \cite{zaher2013optimal}).  As already mentioned,  our approach involves including an $L^\infty$ negative cost term proportional to the maximum/peak inventory levels in the objective function. The difference between the two approaches lies in how the manufacturer's preferences are captured. Constraints impose strict limits (and are apt for existing inventory system), while inclusion of the $L^\infty$-term  allows for a trade-off, leading to potentially better overall outcomes. Such an `optimal policy' is preferable when there is a provision to design the inventory size; of course, this also requires the trade-off factor $\rho$, which reflects the preference of the manufacturer.

We explain the model and methodology in Section  \ref{sec_inventory_model}, then followed by numerical results in Section \ref{sec_numerical_inventory}. Conclusions are in Section~\ref{sec:Conclusions}.

\section{Inventory control with storage design}\label{sec_inventory_model} 
Consider a  manufacturer producing a single product and maintaining an inventory towards the same.  In any inventory system, one needs to control the production rate continuously over the planning time horizon (say $T= [0,t_1]$) to maintain `optimal' inventory level. Further, the demand towards the product is time varying, the customers are sensitive towards the price and hence, the manufacturer also needs to devise a `good' pricing policy. The  higher prices may deter the demand, while  the  lower prices could stimulate it.  Further, the sensitivity towards the price  can also be time varying.

Let $d(t)$ be the demand realized by the manufacturer, by virtue of the   price $p(t)$ quoted by the latter at time $t$. As in \cite{li2015joint}, we consider the following price-demand relation:  
\begin{eqnarray}
     d(t) = \left (\alpha(t) - \beta(t) p(t) \right )^{+}, \ \mbox{ with } z^+ := \max(z, 0) ,
     \label{Eqn_demand}
  \end{eqnarray}
where $\balpha= (\alpha(t), t \in T)$  and $\bbeta = (\beta(t), t \in T)$ represent the maximum possible demand (or the market potential) and the price sensitivity factor respectively. \textit{We assume the manufacturer has the knowledge of $\balpha$ and $\bbeta$, which dictate the future demand curves. But this knowledge may not be  essential---in section \ref{one_control},  some important properties of the optimal policy are derived, by virtue of which, we provide an advice on a `good' (near optimal) policy, which does not require the knowledge of $\balpha$ and $\bbeta$.}

\noindent{\bf Dynamics:}
Let $x(t)$ represent the inventory level and $u(t)$ the production rate at time~$t$. We model the inventory dynamics using the following ODE (with $x(0) = x \ge 0$):
\begin{eqnarray}
    \dot{x}(t)  =  u(t) - (\alpha(t) - \beta(t) p(t))^{+}, 
\label{Eqn_inventory_dynamics}
\end{eqnarray}
  where $\dot{x} = \nicefrac{dx}{dt}$ is derivative of $x$ with respect to $t$ and $x(0)=x$ represents the initial inventory level. In \eqref{Eqn_inventory_dynamics},  $u(t)$ is the rate at which inventory is increasing, and $d(t) = \left (\alpha(t) - \beta(t) p(t) \right )^{+}$ is the demand rate as well as the rate at which the inventory is getting depleted and hence \eqref{Eqn_inventory_dynamics} represents the dynamics.  We basically consider the non-perishable products that does not degrade with time and 
  consider that the demand is served immediately if the inventory is available (i.e, if $x (t) >0$); otherwise, it is backlogged ($x(t) \le 0$) and is served immediately once the  inventory is replenished.
From ODE \eqref{Eqn_inventory_dynamics}, it  may appear that the instantaneous inventory level $x(t)$ does not drive the dynamics. However, the policy tuple  $(\u, \bp)$ typically depends upon $\x$, and this results in feedback to the dynamics.

\noindent{\bf Various costs:} A positive inventory indicates stored products, while the negative levels signify backlogs, indicating unfulfilled demand. When the inventory level $x(t)$ is positive at some time $t$, the manufacturer incurs  an instantaneous holding cost $C_h x(t)^2$ towards the storage; similarly, the backlog instances  result in instantaneous shortage cost $C_s x(t)^2$. The terminal holding and shortage cost factors ($C_h^{(t_1)}, C_s^{(t_1)}$ respectively) can be different from those corresponding to the interim period, possibly due to factors like wastage, repayment, etc. Further, the manufacturer incurs an instantaneous production (or ordering) cost  $au^2(t)$, for producing at rate $u(t)$, where $a$ is a positive constant. 

\noindent{\bf Revenue:} When the demand $d(t)$ is served at time $t$, the manufacturer earns instantaneous profit $p(t)d(t)$. 
 
In all, the revenue generated by the manufacturer over the planning horizon, for any dual policies $(\u, \bp)$,  is given  by  (e.g., see \cite{li2015joint} for such models):
\begin{align}
   & J(0, x; \u, \bp ) 
 =   \int_0^{t_1}\hspace{-2mm}(p(t)(\alpha(t)  
  -\beta(t) p(t))^+ - a u^2(t)  -h(x(t)) ) dt \nonumber\\
& \hspace{50mm}- h_{t_1}( x(t_1) ), \mbox{ with}  \label{Eqn_obj_fun} \\ 
& h(x) := \begin{cases} C_s x^2, \mbox{ if } x <0 \\ C_h x^2, \mbox{ if } x  \ge0,  \end{cases} \hspace{-3mm}
   \mbox{ and } \ 
     h_{t_1} (x) := \begin{cases} C^{(t_1)}_s x^2, \mbox{ if } x  <0 \\ C^{(t_1)}_h x^2, \mbox{ if } x \ge 0.  \end{cases}
  \label{Eqn_terminal_holding_cost}
\end{align}
In the above, ${\bf h}= (h(x))$ and  ${\bf h}_{t_1}= (h_{t_1} (x))$ are the combined inventory holding and shortage cost functions, respectively corresponding to the interim and terminal times.

\noindent{\bf Optimal control:} From \eqref{Eqn_obj_fun}-\eqref{Eqn_terminal_holding_cost},
 it is clear that the inventory optimization problem can be typecast as the following optimal control problem:
 \begin{eqnarray}
     && \sup_{(\u,\bp) \in \cal S }J(0, x; \u, 
    \bp )   \ \   \text{ subject to } \nonumber\\
    &&\dot{x}(t)  =\ u(t) - \left (\alpha(t) - \beta(t) p(t) \right )^{+}, 
\label{Eqn_dyn_without_Linfty}
\end{eqnarray}%
where $\S$ is an appropriate domain and $x(0) = x$ represents the initial condition.

It is obvious  from  \eqref{Eqn_demand}  that the price policy has an upper bound, we require $p(t) \le \nicefrac{\alpha(t)}{\beta(t)}$ for all $t$. We further assume that $u (t)\le \nicefrac {\alpha(t)}{(2 a \beta(t))}$ for all $t$, which is justified by the heuristic arguments as in the footnote\footnote{\label{Footnote_just}Under \textbf{B.1} and for any $\S$ with range of $(u(t), p(t))$ being compact (and sufficiently big) for each $t \in T$,  one can establish that the following 
HJB-PDE  would provide the solution for  \eqref{Eqn_dyn_without_Linfty} (using tools  as in \cite{dhiman2026optimal}, e.g., \cite{fleming2012deterministic, fleming2006controlled}):
$$
\frac{\partial{v}}{\partial{t}} + \sup\limits_{u, p} \left(\frac{\partial{v}}{\partial{x}}\left(u - \left (\alpha(t) - \beta(t) p \right )^{+}\right) +  p(\alpha(t)   -\beta(t) p)^+ - a u^2  -h(x)\right) =0.
$$
Differentiating the term inside supremum wrt $u$ and $p$ and by equating them to zero, the critical points are $p^* = (\frac{1}{2})\frac{\partial {v}}{\partial{x}} + \frac{\alpha(t)}{2 \beta(t)} $ and $u^* = (\frac{1}{2a})\frac{\partial {v}}{\partial{x}}$. Thus $p^* \le \frac{\alpha(t)}{\beta(t)}$ implies  $\frac{\partial {v}}{\partial{x}} \leq \frac{\alpha(t)}{\beta(t)}$, and hence the upper bound on $u$.}.
We thus assume the domain $\S$   is given by:
\begin{eqnarray}
    &&\hspace{-10mm}\S = \bigg \{(\u, \bp) \in {\cal L} (T)\times {\cal L} (T):  0 \leq u(t) \leq \frac{\alpha(t)}{2 a \beta(t)}, \nonumber\\
    && \hspace{20mm}0 \le  p(t) \leq \frac{\alpha(t)}{ \beta(t)} \mbox{ for all } t \in T \bigg\}, 
    \label{Eqn_policy_space}
\end{eqnarray}%
    where ${\cal L} (T)$ is the space of the measurable  functions on $T$ with respect to (wrt) Lebesgue measure, $\muL$. 
We require the following assumption  regarding the   functions describing  the demand~\eqref{Eqn_demand}: \\
\textbf{B.1} \textit{The   functions $\balpha$ and $\bbeta$ are Lipschitz continuous  and   strictly positive  on $T$}.  

Under the above assumption \textbf{B.1}, the policies in domain $\cal S$  are integrable. Further \textbf{B.1} is a typical  assumption considered in the optimal control literature \cite{fleming2006controlled}, which ensures the existence of smooth value function etc., and hence, we proceed with this restriction.


\subsection{Optimal Inventory Storage: $L^{\infty}$ control} \label{sec_Linfty}

We now focus on another important aspect in the inventory management  that of minimizing the peak inventory levels reached during the planning horizon. As already mentioned in the introduction, this aspect is typically neglected in the literature, but represents a very critical design factor, 
 the physical storage capacity---this refers to the maximum amount of inventory that can be stored at any time point in the planned horizon of the operations. Clearly, the physical structure that holds the inventory, once constructed, can not be altered easily. This calls for a careful (or an optimal) design of the storage capacity.

 In literature, one typically solves a constrained optimization problem (objective function as in \eqref{Eqn_obj_fun})  with 
 a hard constraint on the inventory size (e.g., \cite{bitran2003overview, talluri2006theory, zaher2013optimal})---this problem is relevant when the physical storage room is already constructed and hence poses a hard constraint. 
 
 But, we aim to solve a joint optimization problem that optimizes a given weighted combination of the usual inventory costs (as in \eqref{Eqn_obj_fun}) and the physical storage capacity. Such a consideration is well understood to provide a  better policy, but can be implemented only when there is a flexibility to design storage of any required size.
However, interestingly, as illustrated numerically in   subsection \ref{sec_numerical_inventory}, this kind of a flexible design does not advocate exorbitantly large storage---this is true as long as the weightage (or the importance) for the usual inventory cost \eqref{Eqn_obj_fun} is not significantly high.

Towards the above mentioned joint optimization problem, we  consider  ``supremum inventory levels'' or technically introduce the $L^{\infty}$ norm of the inventory levels into the optimization problem \eqref{Eqn_obj_fun}. 
Consider the following parameterized control problem,  with $\rho \ge 0$ quantifying the trade-off between the two  objective functions,  where the control involves tuning two control functions and  with initial condition $ x(0) = x$:
\begin{eqnarray}
    && \hspace{-15mm} v^{(2)}_\rho (0, x) := \sup_{(\u,\bp) \in \cal S } \bigg(J(0, x; \u, \bp )  - \rho \sup_{t \in T} \left\{  x(t) \right\} \bigg)    \mbox{ subject to }   \nonumber\\
  \dot{x}(t) &=&  u(t) - \left (\alpha(t) - \beta(t) p(t) \right )^{+}, \quad  \mbox{ for all } t \in T.  
 \label{Eqn_sup_intro}
\end{eqnarray}
In the above,   $v_\rho^{(2)}(\cdot,\cdot)$ is the value function  (with superscript $.^{(2)}$ indicating  control by two  functions) for the problem that considers $\rho$-weighted combination of the two objectives.  
 The  second objective  $ \sup_{t \in T} \left\{  x(t) \right\} $   embedded into   \eqref{Eqn_sup_intro},   signifies the peak inventory levels achieved over the planned horizon and hence represents the physical storage capacity requirement. In other words, by including this term, we ensure that the control mechanism does not just focus on revenue maximization,  but also    accounts for the costs involved in building the  physical storage. When $\rho = 0$  in  \eqref{Eqn_sup_intro},  the objective function reduces to that in  \eqref{Eqn_obj_fun}, which optimizes the usual inventory costs+profits (also known as, revenue), with zero consideration to the storage capacity. A positive $\rho $  infuses storage capacity cost into the optimal design, whose relative importance   increases with increase in~$\rho$.

 The  combined control problem studied so far in the literature (e.g., \cite{barron1989bellman}) that combines the running cost (like inventory revenue \eqref{Eqn_obj_fun}) and the  $L^{\infty}$ cost  (like peak storage-cost, $\sup_{s} \{ x(s) \}$) has  form   
 
 \vspace{-4mm}
{\small \begin{equation*}
    \sup_{t \in [\tau,t_1]} \left\{ \int_t^{t_1} L(s, x(s), u(s)) ds - \rho   L_\infty (t, x(t),u (t)) \right\} + \Psi(x(t_1)).
     \label{Eqn_avail_Linf}
\end{equation*}}
This    obviously 
cannot capture a more appropriate weighted combination of the two costs   as in   \eqref{Eqn_sup_intro}; such a weighted problem is well-known to capture a part of the Pareto frontier  relevant for muti-objective optimization problems (\cite{ehrgott2005multicriteria}). The smooth variant framework developed in \cite{dhiman2026optimal} will be instrumental in solving \eqref{Eqn_sup_intro}. 

Before proceeding, we  first simplify the problem in \eqref{Eqn_sup_intro} by reducing the domain without losing optimality---the domain $\cal{S}$ is reduced by identifying the relationship between $\bp$ and~$\u$ at optimality.

\subsection{Reduction in Policy Space}
\label{one_control} 
Our problem involves controlling via two  functions simultaneously. The problem simplifies remarkably,   if one can establish some (optimal) relation between the two, which we consider in the immediate next.

\begin{thm}
\label{Thm_one_policy} \textbf{[Reduction in Policy Space]} \textit{Assume \textbf{B.1} and consider  $\S$ as in \eqref{Eqn_policy_space}.\\
(i) Then the problem  \eqref{Eqn_sup_intro} is equivalent to the following control problem, i.e., $v_\rho(0, x) = v^{(2)}_\rho (0, x)$ for all $x$, where: }
\begin{eqnarray}
 v_\rho(0,x) \hspace{-1mm}&\hspace{-1mm}:=\hspace{-1mm}&\hspace{-1mm} \sup_{(\u, \bp) \in {\cal S}}  \bigg(J(0, x; \u, \bp )   - \rho \sup_{t \in T} \left\{  x(t) \right\} \bigg)  \mbox{ subject to } \nonumber \\
  \dot{x}(t) \hspace{-1mm}&\hspace{-1mm}=\hspace{-1mm}&\hspace{-1mm}   u(t) - \left (\alpha(t) - \beta(t) p(t) \right )^{+} \quad   \mbox{ and } \quad  \label{Eqn_ode} \\
  p(t) \hspace{-1mm}&\hspace{-1mm}=\hspace{-1mm}&\hspace{-1mm} a u(t) + \frac{\alpha(t)} {2 \beta(t)}  \quad \forall t \in T.   \quad \label{Eqn_price_prod_relation_cntrl_prblm}
  \end{eqnarray}
\textit{(ii) The problem in part (i) can be rewritten as:}
\begin{eqnarray}
 v_\rho(0,x) \hspace{-1mm}&\hspace{-1mm}:=\hspace{-1mm}&\hspace{-1mm} \sup_{\u \in {\cal S}_{\u }}  \bigg(J(0, x; \u )   - \rho \sup_{t \in T} \left\{  x(t) \right\} \bigg)  \mbox{ subject to } \nonumber \\
  \dot{x}(t) \hspace{-1mm}&\hspace{-1mm}=\hspace{-1mm}&\hspace{-1mm}  u(t)\gamma(t) - \frac{\alpha(t)}{2},  \mbox{ with, }  \nonumber \\
  \gamma(t) \hspace{-1mm}&\hspace{-1mm}=\hspace{-1mm}&\hspace{-1mm}   (1 + a \beta(t)),  \quad \forall t \in T, \label{Eqn_main_reduced}
  \end{eqnarray}
 \textit{ where the  domain }
    \begin{eqnarray} {\cal S}_{\u} 
 := \bigg \{\u  \in {\cal L} (T)   :     u(t)  \in  \left [0, \frac{\alpha(t)}{2 a \beta(t)}\right ], \mbox{ for all } t   \bigg\}, \label{Eqn_u_alone_policy_space} 
  \end{eqnarray}
   and  the objective function modifies  to 
   \begin{eqnarray}
 J(0, x; \u )  &=&  \int_0^{t_1} \bigg( -a u^2(t) \gamma(t)  + \frac{\alpha^2(t)}{4 \beta(t)} -h(x(t))\bigg) d t \nonumber\\
 && \hspace{35mm}- h_{t_1}(x(t_1)). \label{Eqn_only_u_dynamics_objective} 
\end{eqnarray}
\end{thm} 
\noindent {\bf Proof} is in Appendix \ref{sec_app_ch3_Thm_one_policy}. \eop

Thus  the system can now be controlled using only one   control function,  $\u \in \S_\u$ (without loss of optimality).
The   objective function for the reduced problem \eqref{Eqn_main_reduced}-\eqref{Eqn_only_u_dynamics_objective} is a combination of an $L^\infty$ term and an integral term. As discussed in the previous sections, the problem is not a standard optimal control problem. We begin further analysis with existence.

\subsection{Existence} 
\label{subsection_exist_inventory}
We now establish the existence of an optimal solution under relaxed controls for the control problem \eqref{Eqn_main_reduced}-\eqref{Eqn_only_u_dynamics_objective} by verifying that assumptions \textbf{A.0}-\textbf{A.5} of our recent work \cite{dhiman2026optimal} hold under assumption \textbf{B.1}.

To begin with, consider the change of variable, $\tilde{\u}$ defined as  $\tilde{\u} := \nicefrac{ \u  \bbeta  }{\balpha }$ (defined point-wise, recall that by assumption \textbf{B.1}, $\balpha $ and $\bbeta $ are strictly positive). With this transformation,  the problem \eqref{Eqn_main_reduced}-\eqref{Eqn_only_u_dynamics_objective} can be reformulated appropriately. Now,  the control space is given by ${\tilde{\S}}_\u := \left \{\tilde{\u}  : \tilde{u}(t) \in \left [0, \nicefrac{1}{2 a } \right ],  \text{ for all } t   \right\}$, where   $\left [0, \nicefrac{1}{2 a } \right ]$  is compact  and so assumption \textbf{A.0} is satisfied.
Next, let 
$$
f(t,x, {\tilde u})
= {\tilde u} \frac{\alpha(t)}{\beta (t)} \gamma(t) - \frac{\alpha(t)}{2}
$$denote the RHS of the ODE  in \eqref{Eqn_main_reduced}, after using transformed  variable $\tilde{u}$. 
By \textbf{B.1} this RHS   \(f(t,x,\tilde u)\)  is   continuous in $t$, trivially continuous in $\tilde u$  and      is independent of $x$; also observe \(L_{mx}(t,x) \equiv x\) (see \cite[Equation (4)]{dhiman2026optimal} and \eqref{Eqn_main_reduced}).  Thus the  Lipschitz continuity  and uniform boundedness required in assumptions \textbf{A.1} and \textbf{A.3} are trivially satisfied (note here that $\Delta_{mx}  (t,x, {\tilde u}) = f(t,x, {\tilde u}) $, see \cite[Equation (9)]{dhiman2026optimal} ).
Also  \textbf{A.4}  is satisfied as the following is a convex set, for all $t, x$: 
$$
f\left(t, x,  [0, \nicefrac{1}{2a}]\right) =  \left\{ {\tilde u} \frac{\alpha(t)}{\beta (t)} \gamma(t) - \frac{\alpha(t)}{2}  : \tilde{u}\in [0, \nicefrac{1}{2a}]\right\}. 
$$

Since \(\balpha\) and \(\bbeta\) are bounded on the finite interval \(T\), and $\left [0, \nicefrac{1}{2 a } \right ]$ is compact, there exists a constant \(\overline B>0\) such that
$$
\big|f(t,x,\tilde u)\big|\le {\overline B},
\qquad
\forall t\in T,\ \tilde{u} \in  [0, \nicefrac{1}{2a}].
$$
Furthermore, when the initial condition \(x(0)\) is restricted to a neighborhood of zero, or more generally to a compact set, say  \(\mathbb{C}\), the state trajectories has a uniform bound independent of the policy. Indeed, the ODE trajectory under policy $\tilde\u$
\begin{eqnarray*}
|x(t; x(0), \tilde\u)|
&\le& |x(0)|+
\int_0^{t_1}
\big|f(t,x(t),u(t))\big|\,dt \\
&\le& \  \max_{x(0)\in\mathbb C}|x(0)|+{\overline B}t_1, \quad (:= \tilde B) 
\end{eqnarray*}
is bounded uniformly on $T$ and further the bound  is independent of  policy $ \tilde \u \in {\tilde{\S}}_\u$, and so:
$$
\sup\limits_{t\in T, \ x(0)\in\mathbb C, \
\tilde\u\in{\tilde{\S}}_\u}
|x(t;x(0),\tilde\u)|
\le \tilde B<\infty.
$$
Thus when one is only interested in initial conditions coming from a compact set the state trajectories are uniformly bounded irrespective of the policy.

This uniform bound also implies  the boundedness of the functions \(h\) and \(h_{t_1}\) in \eqref{Eqn_terminal_holding_cost}, defining the objective function \eqref{Eqn_only_u_dynamics_objective}. %
 Since these functions are continuous,  their values remain uniformly bounded, again independent of the  policy $\tilde \u$.  Thus we choose to  replace \(h, h_{t_1}\) respectively with
\begin{eqnarray}
    \bar h(x) &:= &\min \{h(x), h_{mx}  \}, \quad h_{mx} := \max_{x \le {\tilde B}} h(x), \mbox{ and } \hspace{5mm}\label{Eqn_defn_barh}\\
    \bar h_{t_1}(x) &:=& \min \{h_{t_1}(x), h_{t_1, mx}  \}, \quad h_{t_1, mx} := \max_{x \le {\tilde B}} h_{t_1}(x) \hspace{5mm}
\label{Eqn_defn_barh_t1}
\end{eqnarray}
to simplify the analysis. 
With this modification, the assumption \textbf{A.5}  is satisfied (in fact \textbf{A.6} is also satisfied).  Assumption \textbf{A.2} is also satisfied as all the involved functions are continuous and $L_{mx} \equiv x$. 

Thus, by  \cite[Theorem 2]{dhiman2026optimal}, an optimal relaxed control exists for the control problem \eqref{Eqn_main_reduced}-\eqref{Eqn_only_u_dynamics_objective}, after replacing $h, h_{t_1}$ with functions in \eqref{Eqn_defn_barh}-\eqref{Eqn_defn_barh_t1}. %
\Forarxiv{In the next section, 
we derive some structural  properties of the optimal policies under some additional conditions related to existence. These additional conditions are required only for the next section.

 \subsection{Qualitative analysis of optimal policy}

Under the  following additional assumption, we obtain some interesting properties of the optimal policy: \\
\textbf{B.2} \textit{Assume  the existence of the optimal policy  $\u^*$  among the pure controls, which is further  in the interior of $\S_\u$:  there exists an $\eta > 0$ such that  } $
 u^*(t) \in \big[\eta, \nicefrac{\alpha(t)}{2 a \beta(t)} - \eta \big] 
 $, \textit{for each} $t \in T$.   
\begin{thm}\label{thm_w_dec_inc}\textbf{[Characteristics of Optimal Policy]} 
        Assume \textbf{B.1}-\textbf{B.2}.  Let  $\x^*$ be the state trajectory corresponding to the optimal policy $\u^*$. Define  
        $$
        w(t, t') := \int_{t}^{t'} u^*(s) ds,  
         \quad  \mbox{ for all } t, t' \in T.
        $$
(i) If at some time $t \in (0, t_1)$, 
            \begin{itemize}
                \item[(a)] $x^*(t) >0$, then there   exists a $\bar \delta > 0$  such that 
                 $
                w(t - \delta, t) \le w(t , t + \delta) \mbox{ for all }   \delta \le \bar \delta, 
                $
                \item[(b)] $x^*(t) <0$, then there exists a $\bar \delta> 0$  such that 
                $
                w(t - \delta, t) \ge w(t , t + \delta) \mbox{ for all }   \delta \le \bar \delta.
                $
            \end{itemize}
(ii)  There exists a $\delta >0$, such that $x^{*}(t)\le 0$,  for all $t \in [t_1 - \delta,  t_1]$.
 \end{thm}
 \noindent{\bf Proof} is available in Appendix \ref{sec_app_ch3_ch_opt_policy}.    \eop

From part $(i)$ of Theorem \ref{thm_w_dec_inc}, the optimal production rate locally increases (or decreases) with time when the inventory level at that instant is strictly positive (or negative respectively). 
This theorem suggests a procedure for manufacturers: 
when the inventory has surplus quantity, it is beneficial  to choose the production rate conservatively: beneficial  to choose a smaller rate initially and an increased rate at a later instant.  On the other hand, when the manufacturer is facing a shortage (negative inventory), it is advantageous to consider a higher production rate initially and decrease later. \textit{This kind of a strategy advice is probably beneficial even in the scenarios with partial/zero information on future demand trajectories $\balpha$ and $\bbeta$.} We discuss  more on this  at  end of subsection~\ref{sec_numerical_inventory}, which provides numerical analysis.

From part $(ii)$ of Theorem \ref{thm_w_dec_inc}, 
   it is  advantageous  to have shortage towards the end of the horizon, this  is true for any set of parameters. In other words, even when the shortage cost parameters are `significantly higher' than  the  holding cost parameters and~$\rho$, the optimal policy still ensures shortage towards the end (probably the amount of shortage $x^*(t_1)$ towards the end is smaller, for  larger values of $C_s, C_s^{(t_1)}$). This   implies another strategy advice: \textit{it is beneficial to ensure zero inventory towards the end and repay any unfulfilled orders}. }

 We now derive probably an  obvious, yet important result in the following: 

\begin{propsn} Assume \textbf{B.1}.  Define $y_\rho^*(t_1) := \sup_{t \in T} \left\{  x_\rho^*(t) \right\}$ at the optimal state trajectory  $\x_\rho^*$. Then $y_\rho^*(t_1)$ decreases with increase in $\rho$.
\label{thm_y_terminal_dec_with_sigma} 
\end{propsn} 
 \noindent{\bf Proof} is available in Appendix \ref{sec_app_ch3_thm_y_terminal_dec_with_sigma}.    \eop
 
Observe that $y_\rho^*(t_1)$ defined in the above proposition represents the optimal storage capacity requirement for the given trade-off factor $\rho.$ As proved above, with increased $\rho$, one can have an optimal policy that requires a smaller physical storage capacity. However, the immediate question is regarding the price one has to pay in terms of the revenue generated, in exchange for smaller physical storage capacity. Our next quest is to analyze  the effect of the $L^{\infty}$ term on the revenue in~\eqref{Eqn_main_reduced}. To be more precise, our goal is to explore how varying the importance attributed to the $L^{\infty}$ cost (via $\rho$) affects the `optimal' revenue component $J^* := J(0, x; \u^*)$ (defined using optimal $ \u^*$ trajectory)  in \eqref{Eqn_main_reduced}. Interestingly, we will discover that despite the increase in $\rho$, one would not sacrifice much in terms of the revenue $J^* $, while the physical storage capacity significantly reduces.  Towards illustrating the same and more, we next proceed with a  smooth   variant of \eqref{Eqn_main_reduced}, obtained using the framework of  \cite[section IV]{dhiman2026optimal}.

\subsection{Conversion to standard problem}
\label{subsec_inventory_smooth}
It is clear that the inventory control problem \eqref{Eqn_main_reduced} is an example of the general problem considered   in \cite[Equation 4]{dhiman2026optimal}. 
The corresponding smooth variant  can be 
obtained by  following the smooth approximation procedure described in  \cite[Equation 21]{dhiman2026optimal}, where function \(\psi_\delta\) is defined using any Lipschitz continuous increasing function \(\phi:[0,1]\to[0,1]\) satisfying \(\phi(0)=0\) and \(\phi(1)=1\) (see again \cite[Equation 21]{dhiman2026optimal}). In particular,  we choose \(\phi\) as an identity function, and  obtain the following simplified smooth control problem:
\begin{align}
  v_\rho(0,x, y) &:= \sup_{\u \in {\cal S}_\u }  \bigg(J(0, x; \u )   - \rho   y(t_1)   \bigg) \quad  \mbox{ subject to 
   }  \nonumber\\ 
 \dot{x}(t)  &= u(t) \gamma(t) -  \frac{\alpha(t)}{2},  \nonumber\\
   \dot{y}(t) &=  \left (u(t) \gamma(t) -  \frac{\alpha(t)}{2} \right )^+  \psi_\delta ( x(t)- y(t)),  \quad  \mbox{ where } \nonumber 
\\    \gamma(t) &:= 1 + a \beta(t), \quad  \mbox{ for all }t \in T \mbox{ and }\nonumber  
\\
J(0, x; \u ) &= \int_0^{t_1} \bigg( -a u^2(t) \gamma(t)  + \frac{\alpha^2(t)}{4 \beta(t)} 
    -\bar h(x(t))\bigg) dt \nonumber \\
    & \hspace{40mm}- \bar h_{t_1}(x(t_1)),  
\label{Eqn_final_control_problem}
\end{align}%
where $\bar h(\cdot)$ and $\bar h_{t_1} (\cdot)$ are defined in \eqref{Eqn_defn_barh} and \eqref{Eqn_defn_barh_t1}. 

In Subsection \ref{subsection_exist_inventory}, we established that assumptions \textbf{A.0}--\textbf{A.5} hold under \textbf{B.1}. 
Consequently, by \cite[Lemma 2]{dhiman2026optimal}, the smooth problem \eqref{Eqn_final_control_problem} admits an optimal solution within the class of pure controls $\S_\u$ (noting that the transformation $\u \mapsto \tilde \u$ is bijective). 
Furthermore, as previously noted, assumption \textbf{A.6} of  \cite{dhiman2026optimal} is also satisfied. 
Therefore,  the value function of the smooth problem \eqref{Eqn_final_control_problem} is the unique viscosity solution to the following HJB PDE:
\begin{eqnarray} 
 &&\hspace{-5mm}\frac{\partial v}{\partial t}  + \sup_{u \in \left [0, \frac{\alpha}{2 a \beta}  \right ]}  H\left(t, x, y, u, \frac{\partial v}{\partial x} ,\frac{\partial v}{\partial y} \right) = 0, \mbox{ where } \nonumber  \\ 
 &&\hspace{-5mm} H(t, x, y, u, p, q) := p\dot{x} + q \dot{y} + \frac{\alpha^2(t)}{4 \beta(t)} -a u^2 \gamma(t)  - \bar h(x), \mbox{ with } \nonumber \\
  &&\hspace{-5mm}v(t_1, x, y) = - \rho y - \bar h_{t_1}(x). \label{Eqn_hjb}
  \end{eqnarray}

\vspace{-2mm}
\section{Numerical Results} 
\label{sec_numerical_inventory}

We have established property   of the direct problem \eqref{Eqn_main_reduced}  in Theorems \ref{Thm_one_policy}. We now turn to the numerical analysis to gain further insights into the problem.

We use the Pontryagin Maximum Principle (see e.g., \cite{fleming2006controlled})  to numerically solve the problem. Using  the Hamiltonian $H$ in \eqref{Eqn_hjb},   one has to solve   the following  co-state equations and  compute the solution $(\lambda_x (\cdot), \lambda_y(\cdot), x(\cdot), y(\cdot))$, see \eqref{Eqn_final_control_problem}: 
\begin{eqnarray}
    \frac{d \lambda_x}{d t} &=& 2 C_s x \mathds{1}_{\{x<0\}}+2 C_h x \mathds{1}_{\{x>0\}}   - \frac{\lambda_y}{\delta} \dot{x}\mathds{1}_{\{ x- y \in [-\delta, 0]\}}, \nonumber\\ 
\frac{d \lambda_y}{d t} &=&   \frac{\lambda_y}{\delta} \dot{x} \mathds{1}_{\{ x- y \in [-\delta, 0]\}}, 
\label{Eqn_pontygn_principle}
\end{eqnarray}
with   terminal conditions:
{\small$
\lambda_x(t_1) = -C_s^{(t_1)} x \mathds{1}_{\{x<0\}}-C_h^{(t_1)} x \mathds{1}_{\{x > 0\}}, $ $
\lambda_y(t_1)  = - \rho,
$}
and then the optimal policy is:
$$ 
u^* =\max\bigg\{0,  \min\bigg\{\frac{\alpha}{2 a \beta}, \frac{\lambda_x}{2a} + \frac{\lambda_y}{2a} \mathds{1}_{\{u^* \gamma \ge \frac{\alpha}{2}, \  x- y \in [-\delta, 0]\}} \bigg\}\bigg\}.
$$

We set $\delta = 0.01$ in all examples, except for some tables; here we also study the effect of smooth approximation (along with others), by comparing across different $\delta$. First  note that the iterative procedure  converges in all the case studies.  
We begin with the study of the  system at different  trade-off factors $\rho$.

\subsubsection{Effect of $L^{\infty}$ cost, as $\rho$ varies} 
\label{subsec_effect_of_infty}
\begin{figure*}[h]
\hspace{-20 mm}
    \centering
    \begin{minipage}{5cm}
\includegraphics[trim = {4.5cm 1.8cm 2cm 2.7cm}, clip, width = 6cm, height = 3.8cm]{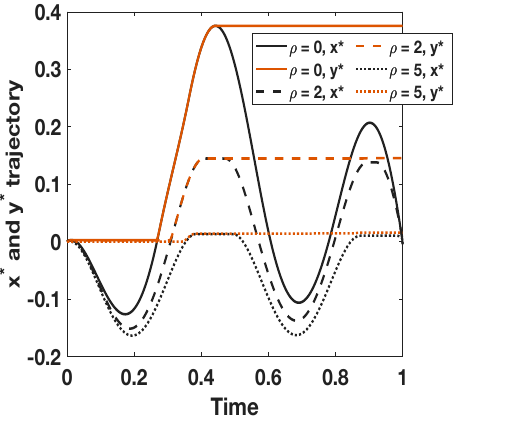}
    \end{minipage}
    \hspace{0.1cm}
        \begin{minipage}{5cm}
\includegraphics[trim = {0cm 0cm 0cm 0cm}, clip, width = 4.6cm, height = 3.7cm]{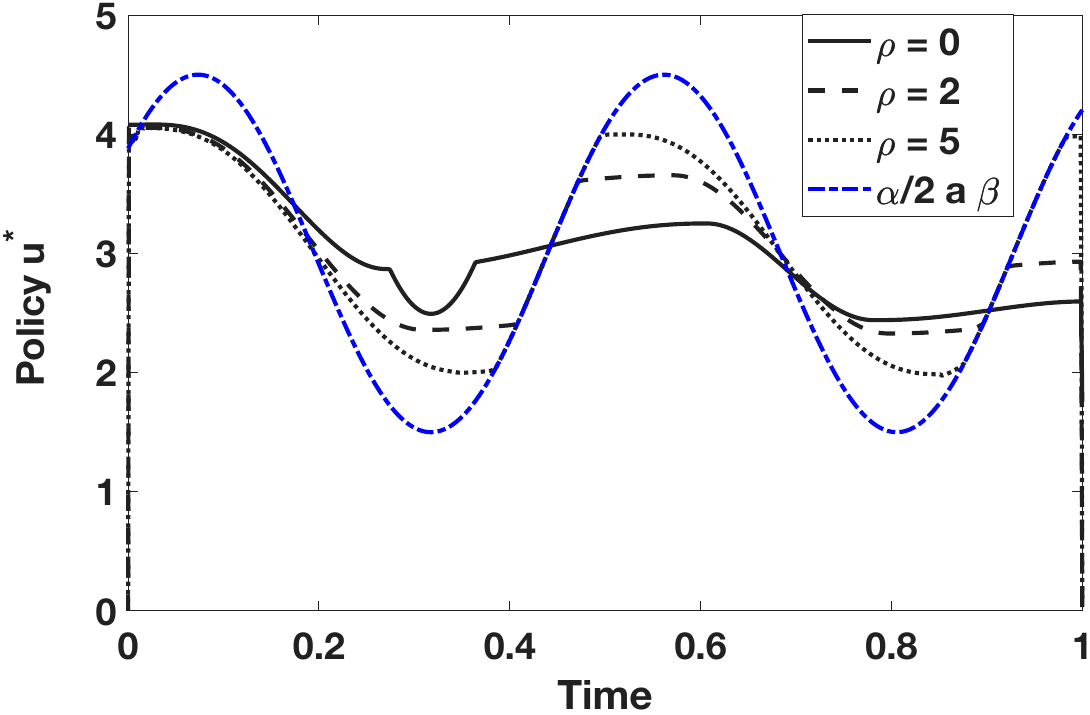}
    \end{minipage}
    \hspace{0.1cm}
        \begin{minipage}{5cm}

{\small\begin{tabular}{|l|l|l|l|l|}
\hline   \multicolumn{2}{|l|}{} & \multicolumn{1}{c|}{$\delta=0.5$} &{$\delta=0.05$} & $\delta = 0.01$ 
\\ \hline
\multirow{2}{*}{$\rho =0$}  & $J^*$   & 11.391 & 11.391   & 11.3888  
\\ \cline{2-5}   & $y^*(T)$   &0.4824  & 0.3737    & 0.3646     
\\ \hline
\multirow{2}{*}{$\rho = 2$} & $J^*$  & 10.5402 & 10.8986    & 10.9206 
\\ \cline{2-5}   & $y^*(T)$  & 0.2933 & 0.16319  & 0.1445 
\\ \hline
\multirow{2}{*}{$\rho =5$}  & $J^*$  & 10.2478  & 10.6208   & 10.7189  
\\ \cline{2-5}    & $y^*(T)$ &  0.0358 & 0.0256  & 0.01564  
\\ \hline
\end{tabular}}
    \end{minipage}
    \vspace{0mm}
    \caption{ 
Case-study 1: State trajectories $(x^*, y^*)$, optimal policy $u^*$  and optimal $(J^*, y^*)$  as~$\rho$  varies. 
}
    \label{Fig_c1_com_state_policies}
\end{figure*}

\begin{figure*}[h]
\vspace{1mm}
\hspace{-15mm}
    \centering
    \begin{minipage}{5cm}
\includegraphics[trim = {0cm 0cm 0cm 0cm}, clip, width = 4.7cm, height = 4cm]{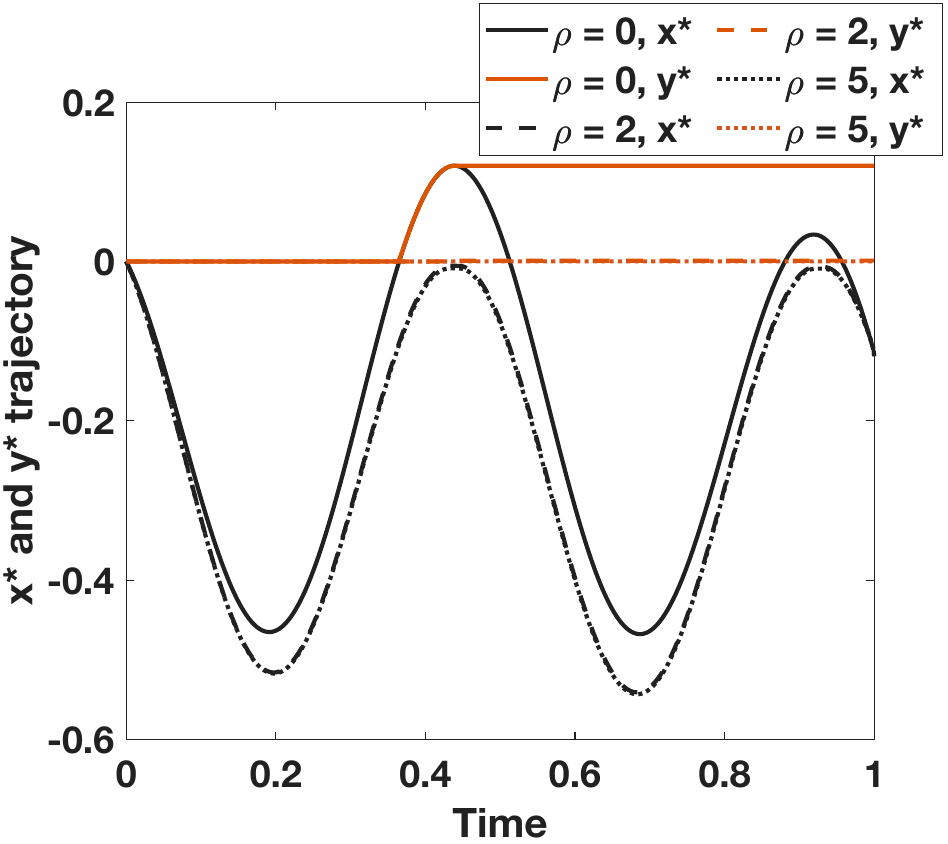}
    \end{minipage}
    \hspace{0.0cm}
        \begin{minipage}{5cm}
\includegraphics[trim = {0cm 0cm 0cm 0cm}, clip, width = 4.6cm, height = 3.4cm]{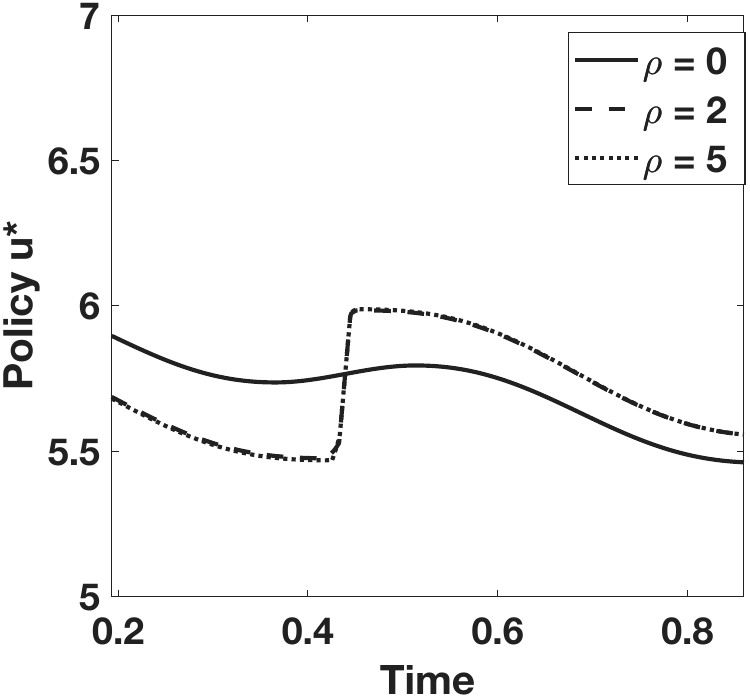}
    \end{minipage}
    \hspace{0.1cm}
        \begin{minipage}{5cm}
{\small \begin{tabular}{|l|l|l|l|l|}
\hline   \multicolumn{2}{|l|}{} & \multicolumn{1}{c|}{$\delta=0.5$} & $\delta = 0.05$ & $\delta = 0.01$ 
\\ \hline
\multirow{2}{*}{$\rho =0$}  & $J^*$   & 102.129& 102.130   & 102.130  
\\ \cline{2-5}   & $y^*(T)$  & 0.1258 & 0.1247    & 0.120     
\\ \hline
\multirow{2}{*}{$\rho = 2$} & $J^*$  & 102.040 & 102.077 &102.0847     
\\ \cline{2-5}   & $y^*(T)$   & 0.0114  & 0.001 & 0.00005
\\ \hline
\multirow{2}{*}{$\rho =5$}  & $J^*$   &  102.026   & 102.076 & 102.084 
\\ \cline{2-5}    & $y^*(T)$ & 0.0015  & 0.00001 & 0.0001  
\\ \hline
\end{tabular}}
    \end{minipage}
    \vspace{-4mm}
    \caption{ 
Case-study 2 with a more relaxed market. }
    \label{Fig_c2_com_state_policies}
\end{figure*}

We then study the impact of $\rho$ on the solutions and the optimal revenues 
in two varieties of markets distinguished via characteristics like the sensitivity parameters and the shortage-cost parameters. We basically consider two extreme sets of parameters to cover interesting varieties of examples. In all the examples,  a wide range of $\rho$ starting from $0$ is considered to cover all possible trade-offs (i.e., to study problems with different levels of importance to inventory size). \\
 \textbf{Case Study 1:} In the first case study, we consider a market that is highly sensitive 
 and where the shortage costs are extremely high. Such a case study is representative of a `\textit{tough}' market, where 
we consider specific parameters: $x(0)=y(0)= 0, C_h = 3, C_s = 40, a = 0.6, t_1 = 1, \beta(t) = 2.5, \alpha(t) = 15 + 4.5 \sin(0.2 \pi + 4.1 \pi t), C_s^{(t_1)} = 410$ and $ C_h^{(t_1)} = 6 $ (see Figure \ref{Fig_c1_com_state_policies}). We use short notation $J^*$ 
in the  tables to represent the value  $J^* = J(0, 0;u^*)= v( 0, (0,0))+ \rho y^*(t_1)$  of  \eqref{Eqn_final_control_problem}.\\
 \textbf{Case Study 2:} In contrast, this scenario represents a more relaxed market environment, with lower shortage costs and sensitivity parameters:
 $x(0)=y(0)= 0, C_h = 3, C_s = 2, a = 0.6, t_1 = 1, \beta(t) = 0.5, \alpha(t) = 15 + 4.5 \sin(0.2 \pi + 4.1 \pi t), C_s^{(t_1)} = 30$ and $ C_h^{(t_1)} = 6 $ (see Figure~\ref{Fig_c2_com_state_policies}).

In Figure \ref{Fig_c1_com_state_policies}-\ref{Fig_c2_com_state_policies}, we solve the problem \eqref{Eqn_final_control_problem} using Pontryagin Maximal Principle \eqref{Eqn_pontygn_principle}. Without including $L^{\infty}$ cost, i.e., when $\rho = 0$, we observe from both the figures that the revenue  is the maximum (see the tables included in the figure environment). Further,  the maximum inventory level or the $L^\infty$ cost  is the largest  (the black continuous curve represent the optimal instantaneous inventory levels ${\bm x}^*$, while the red curve depict ${\bm y}^*$, see also the tables).  
In the left sub-figures
  of  Figures \ref{Fig_c1_com_state_policies}-\ref{Fig_c2_com_state_policies}, the plots depict  ${\bm x}^*$ and ${\bm y}^*$ trajectories for different values of $\rho$. The inventory ${\bm x}^*$ exhibits significant variations with zero  importance to $L^\infty$ cost  (i.e., with $\rho = 0$), indicating high volatility in the inventory trajectory. However, as $\rho$ increases, the variations in the inventory  decrease (dash and dotted lines).

 \begin{wrapfigure}{r}{0.22\textwidth}
\begin{center}
\vspace{-7mm}
\includegraphics[trim = {0cm 0cm 0cm 0cm}, clip, width=3.8cm, height=3.2cm]{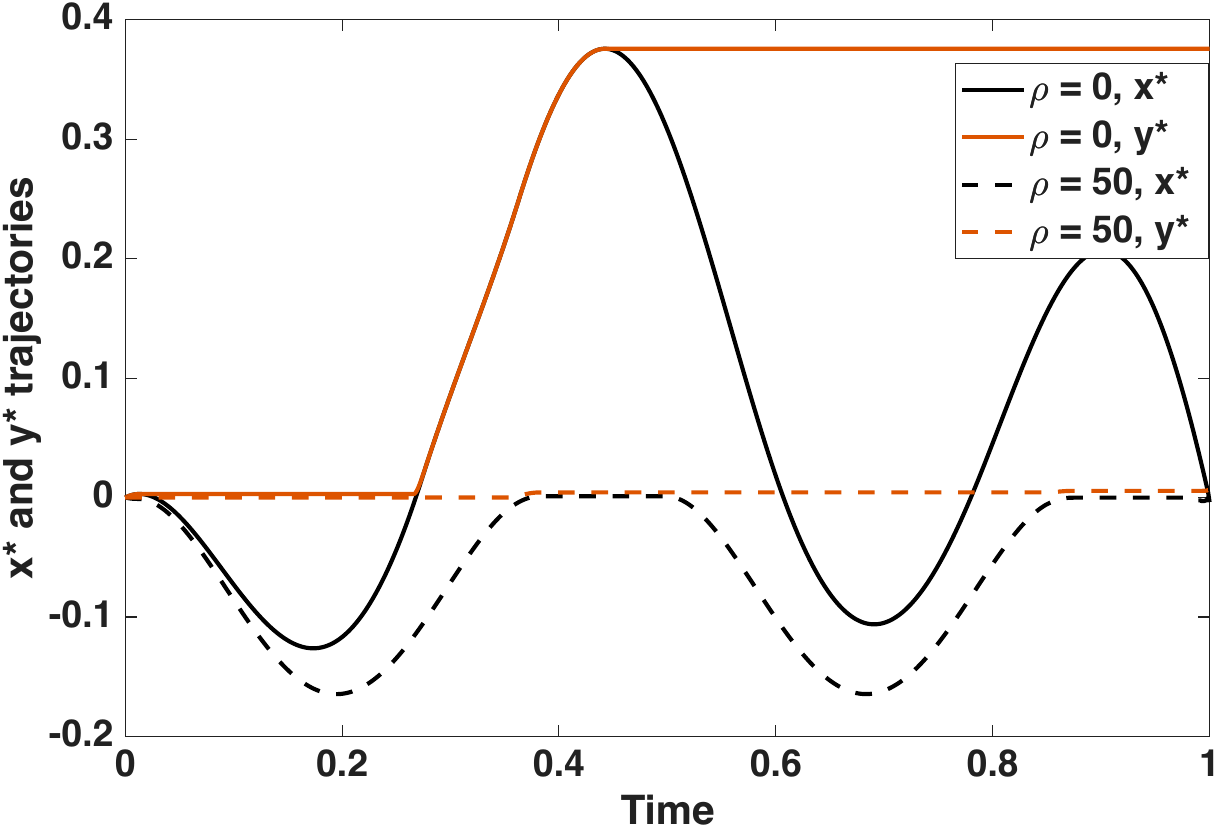}
 \vspace{-2mm}
\caption{$x^*, y^*$ trajectories}
\vspace{-5mm}
\end{center}
\label{Fig_comp_traj}
\end{wrapfigure}
The right sub-figures of Figures \ref{Fig_c1_com_state_policies}-\ref{Fig_c2_com_state_policies},   display    the  optimal policy (black curves). 
We also display the  $\balpha/2a\bbeta$ trajectory for comparison purposes.  When $x^*(t)$ is negative,  the production rate $u^*(t)$ locally decreases, and increases locally when $x^*(t)$ is positive;  as seen from the figures, this is true at majority of time points $\{t\}$, except for the ones at the switching points. 
\Forarxiv{Thus the smooth approximate solution also satisfies Theorem \ref{thm_w_dec_inc}  majorly.} 

 As $\rho$ increases, we observe a trend: without compromising on revenue by more than $6\%$ in case-study 1 (and negligible for case-study 2 with less sensitive market), the 
 maximum inventory-level $(y^*(t_1))$ reduces significantly.
 In fact, at $\rho = 5$, $y^*(t_1)$ is almost zero. 
Further,  the   fluctuations in the optimal inventory trajectory reduce with increase in $\rho$. This trend becomes more prominent  with higher values of $\rho$. In Figure 3, with $\rho = 50$, we in fact notice negligible fluctuations. The  tables in both the figures \ref{Fig_c1_com_state_policies} and \ref{Fig_c2_com_state_policies} also show the results with different $\delta$, we observe various quantities are converging as $\delta$ decreases, indicating consistent results with smooth approximation.
\Forarxiv{This further strengthens the strategy advice provided after Theorem \ref{thm_w_dec_inc}---the manufacturer should maintain inventory as close to zero as possible throughout the planning horizon. }

\section{Conclusions}\label{sec:Conclusions} The maximum-level reached by the inventory during the planning horizon is representative of the storage capacity requirement and the focus of this paper is to include the same in the inventory optimization problem. This addition leads to the  non-standard control problem, of the form introduced in our recent paper \cite{dhiman2026optimal}, involving a combination of classical integral term and an $L^\infty$
 term. 
 
We first  establish a
relationship between the optimal pricing and production-rate
policies, which reduced the dimension of the problem. Using the existence results established in  \cite{dhiman2026optimal}, we obtain an optimal relaxed control and, through a smooth approximation procedure, construct an $\epsilon$-optimal pure control for the reduced problem.
 \Forarxiv{We also  illustrate  certain salient features of the optimal policy, under some additional conditions. At the instants of positive inventory the optimal production rate should be increasing locally in time and locally decreasing with shortage. The optimal inventory will not end with excess.} 

By numerically solving the smooth approximate problem, and  obtain that with increased importance to 
$L^\infty$
 term or storage capacity cost, the fluctuations in optimal inventory trajectory reduce (and the reduction in revenue is not significant, around 
6\%
 in highly sensitive markets and negligible otherwise). The overall directive of this paper is: the optimal policy for the manufacturer is to manage the inventory near zero. And such a policy can be implemented with or without information about the future demand curve.

Overall, the paper deals with a non-standard control problem for inventory control. By introducing an auxiliary variable,   the non-standard problem convert to the classical framework. However, the resultant has discontinuities that require significant technical attention. It would be interesting to solve such problems in general, and probably can identify if there exists a special dynamic programming (Bellman) equation for such problems.

\bibliographystyle{IEEEtran}
\bibliography{references}

\begin{thebibliography}{10}
\providecommand{\url}[1]{#1}
\csname url@samestyle\endcsname
\providecommand{\newblock}{\relax}
\providecommand{\bibinfo}[2]{#2}
\providecommand{\BIBentrySTDinterwordspacing}{\spaceskip=0pt\relax}
\providecommand{\BIBentryALTinterwordstretchfactor}{4}
\providecommand{\BIBentryALTinterwordspacing}{\spaceskip=\fontdimen2\font plus
\BIBentryALTinterwordstretchfactor\fontdimen3\font minus \fontdimen4\font\relax}
\providecommand{\BIBforeignlanguage}[2]{{%
\expandafter\ifx\csname l@#1\endcsname\relax
\typeout{** WARNING: IEEEtran.bst: No hyphenation pattern has been}%
\typeout{** loaded for the language `#1'. Using the pattern for}%
\typeout{** the default language instead.}%
\else
\language=\csname l@#1\endcsname
\fi
#2}}
\providecommand{\BIBdecl}{\relax}
\BIBdecl

\bibitem{chen2004coordinating}
X.~Chen and D.~Simchi-Levi, ``Coordinating inventory control and pricing strategies with random demand and fixed ordering cost: The infinite horizon case,'' \emph{Mathematics of operations Research}, 2004.

\bibitem{chen2012pricing}
------, ``Pricing and inventory management,'' 2012.

\bibitem{elmaghraby2003dynamic}
W.~Elmaghraby and P.~Keskinocak, ``Dynamic pricing in the presence of inventory considerations: Research overview, current practices, and future directions,'' \emph{Management science}, vol.~49, no.~10, 2003.

\bibitem{whitin1955inventory}
T.~M. Whitin, ``Inventory control and price theory,'' \emph{Management science}.

\bibitem{tu2013introductory}
P.~N. Tu, \emph{Introductory Optimization Dynamics: optimal control with economics and management science applications}.\hskip 1em plus 0.5em minus 0.4em\relax Springer Science \& Business Media, 2013.

\bibitem{zaher2013optimal}
H.~Zaher, N.~R. Said, and T.~T. Zaki, ``Optimal control theory approach to solve constrained production and inventory system,'' \emph{International Journal}, vol.~3, no.~12, 2013.

\bibitem{polak1993method}
E.~Polak, T.~H. Yang, and D.~Q. Mayne, ``A method of centers based on barrier functions for solving optimal control problems with continuum state and control constraints,'' \emph{SIAM Journal on Control and Optimization}, vol.~31, no.~1, pp. 159--179, 1993.

\bibitem{schiela2009barrier}
A.~Schiela, ``Barrier methods for optimal control problems with state constraints,'' \emph{SIAM Journal on Optimization}, vol.~20, no.~2, pp. 1002--1031, 2009.

\bibitem{li2015joint}
S.~Li, J.~Zhang, and W.~Tang, ``Joint dynamic pricing and inventory control policy for a stochastic inventory system with perishable products,'' \emph{International Journal of Production Research}, 2015.

\bibitem{barron1989bellman}
E.~Barron and H.~Ishii, ``The {Bellman} equation for minimizing the maximum cost.'' \emph{Nonlinear Analysis: Theory, Methods \& Applications}.

\bibitem{barron1999viscosity}
E.~Barron, ``Viscosity solutions and analysis in $l^\infty$,'' \emph{Nonlinear analysis, differential equations and control}, pp. 1--60, 1999.

\bibitem{dhiman2026optimal}
M.~Dhiman, V.~Kavitha, and N.~Hemachandra, ``Optimal control with ${L}^\infty$ and integral cost functionals,'' \emph{arXiv preprint arXiv:2608.14316}, 2026.

\bibitem{fleming2012deterministic}
W.~H. Fleming and R.~W. Rishel, \emph{Deterministic and stochastic optimal control}.\hskip 1em plus 0.5em minus 0.4em\relax Springer Science \& Business Media, 2012, vol.~1.

\bibitem{fleming2006controlled}
W.~H. Fleming and H.~M. Soner, \emph{Controlled Markov processes and viscosity solutions}.\hskip 1em plus 0.5em minus 0.4em\relax Springer Science \& Business Media, 2006.

\bibitem{bitran2003overview}
G.~Bitran and R.~Caldentey, ``An overview of pricing models for revenue management,'' \emph{Manufacturing \& Service Operations Management}.

\bibitem{talluri2006theory}
K.~T. Talluri and G.~J. Van~Ryzin, \emph{The theory and practice of revenue management}.\hskip 1em plus 0.5em minus 0.4em\relax Springer Science \& Business Media, 2006, vol.~68.

\bibitem{ehrgott2005multicriteria}
M.~Ehrgott, \emph{Multicriteria optimization}.\hskip 1em plus 0.5em minus 0.4em\relax Springer Science \& Business Media, 2005.

\end{thebibliography}

\section{Appendix}
\label{sec_app_ch3}
\section*{Proof of Theorem \ref{Thm_one_policy}} 
\label{sec_app_ch3_Thm_one_policy}
For any given policy $(\u, \bp)$, 
define 
\begin{equation*}
    z(t):= u(t) -  (\alpha(t) - \beta(t) p(t))^+ \mbox{, for all } t \in T.  
\end{equation*}
Observe from  the above,  many pairs of $(\u, \bp)$ can lead to same $\z = (z(s), s\in T )$ and  each such pair leads to same state trajectory, call it  $\x_\z $,   as seen from \eqref{Eqn_ode}. In other words,  the evolution of the state trajectory does not change with particular choice of the pair $(\u, \bp)$, as long as the pair lead to the same $\z$ (we represent this briefly as $(\u, \bp) \leadsto \z$).  

Let $\S_\z := \{ (\u, \bp) \in \S : (\u, \bp) \leadsto \z\}$, for any consolidated policy   trajectory $\z$,  represent the subset of $\S$ that leads to  $\z$ and hence to same state trajectory $\x_\z$.
We will now choose   the optimal pair  $(\u^*_\z, \bp^*_\z)$ among each sub-class   $  \S_\z$,  i.e.,  for each $\z$, which optimizes or satisfies: 
$$
J(0, x; \u^*_\z, \bp^*_\z) = \sup_{ (\u, \bp) \in \S_\z}  J(0, x; \u, \bp).  
$$
We will show below that such an sub-optimal  pair    satisfies the relation given in \eqref{Eqn_price_prod_relation_cntrl_prblm} for each $\z$, which completes the proof of part $(i)$. Fix one such $\z.$

From the structure of     the cost functional \eqref{Eqn_sup_intro},  the only terms that depend upon the specific details of $(\u, \bp)$ are instantaneous terms, which at $t$ are given by:  $O_t (\u,\bp) :=  p(t)(\alpha(t)   -\beta(t) p(t))^+ - a u^2(t).$ The rest of the terms that define the objective functional (including the $L^\infty$-cost) in \eqref{Eqn_sup_intro} depend only upon state trajectory $\x_\z$.  
Thus it 
suffices to solve the following optimization problem independently for each $t \in T$ and show that the optimal pair $(\u_\z^*, \bp_\z^*)$ satisfies the relation in \eqref{Eqn_price_prod_relation_cntrl_prblm} at the same $t$  (see \eqref{Eqn_policy_space}): 
\begin{eqnarray}
&& \hspace{-44mm}(u_\z^*(t), p_\z^*(t)) \in   \arg \max_{ u \in [0, \bar u(t)], p\in [0, \bar p(t)] }  O_t ( u,p) \ \mbox{ subject to }\label{Eqn_Ot_max}\\
u  -  (\alpha(t) - \beta(t) p )^+ &=& z(t)
\quad \mbox{ with } \nonumber\\
O_t (u,p ) &:=&  p(\alpha(t)   -\beta(t) p)^+ - a u^2, \nonumber\\
\bar u(t) &=& \frac{\alpha(t)}{2 a \beta (t)},  \mbox{ and }  \bar p(t) = \frac{\alpha(t)}{\beta (t)}.  \nonumber  
\end{eqnarray}
The above optimization is  solved in two sub-domains: \\
{\bf Case 1:} If $p = \bar p(t)$,  then $u = z(t)$. For such $(u,p)$, 
$$
O_t (u,p )   = - a z(t)^2.
$$
{\bf Case 2:} If $p < \bar p(t)$, then $z(t) = u - (\alpha(t) - \beta (t) p)$ and then $u = z(t) + (\alpha(t) - \beta (t) p)$; thus  in this case: 
$$
O_t(u,p ) = p \big(\alpha(t) - \beta (t) p\big) - a \big(z(t) + (\alpha(t) - \beta (t) p)\big)^2.  
$$
Towards finding the optimal in this sub-domain, first observe:
$$
\frac{\partial O_t (u,p)}{\partial p } = 
\alpha(t) - 2 \beta (t) p + 2 a \beta (t) \big(z(t) + (\alpha (t) - \beta (t) p \big). 
$$
Now
obtain the zero of the above derivative and then 
 by concavity we have (and as the derivative is positive at $p = 0$). 
\begin{eqnarray}
    p^*_\z (t)
    &=&  \min \left \{ \frac{ \alpha (t)  + 2 a \beta (t) (z(t) + \alpha (t) ) } {2 \beta (t) + 2 a \beta(t)^2}, \frac{\alpha(t)}{\beta (t) } \right \} \nonumber\\
    &=&   \frac{ \alpha (t)}{\beta(t) }  -\frac{  ( \alpha(t)  - 2 a \beta (t) z(t)   )^+  } {2 \beta (t) + 2 a \beta(t)^2}.    
    \label{Eqn_p_star}
\end{eqnarray}
By concavity again, the optimal over the entire domain (by considering Case 1 and 2) is given by $ p^*_\z (t) $ as above and 
\begin{eqnarray}
u^*_\z (t) &=& z(t) +  \mathds{1}_{\{p^*_\z (t) < \bar p(t)\}}  (\alpha(t) - \beta (t) p^*_\z (t) ) \nonumber\\
&=&  z(t) +  \mathds{1}_{\{p^*_\z (t) < \bar p(t)\}}  \frac{  \alpha(t)  - 2 a z(t) \beta(t) }{2 + 2 a \beta (t) } \nonumber\\
&=&   z(t) +     \frac{  (\alpha(t)  - 2 a z(t)  \beta(t))^+   }{2 + 2 a \beta (t) }. 
\label{Eqn_u_star}
\end{eqnarray}
The above pair would represent the actual optimizers of \eqref{Eqn_Ot_max}, if they additionally satisfy the constraints. 

First consider the case when 
    $z(t) < \alpha(t)/ (2 a \beta (t) )$, then  using \eqref{Eqn_p_star} and  \eqref{Eqn_u_star}, 
we clearly  have the following, which reflects the relation in \eqref{Eqn_price_prod_relation_cntrl_prblm} and the required constraint on~$u_\z^*(t)$:  
\begin{equation*}
p^*_\z (t) = a u^*_\z (t) + \frac{\alpha(t)} {2 \beta(t)} \mbox{ and further that } u^*_\z (t)  <  \frac{\alpha(t)} { 2 a \beta (t)} = \bar u_t. 
\end{equation*}
In the other case,  we have: $u_\z^*(t) = z(t)$ and $p_\z^*(t) =   \bar p(t)$; in this case by definition of~$\S$, the only allowed value of    $z(t) = \bar u(t)$ and then again $(u_\z^*(t) , p_\z^*(t) )$
 satisfy \eqref{Eqn_price_prod_relation_cntrl_prblm}.
Thus, part~$(i)$ is   proved.

Towards part (ii), define the following mapping $\p:  \S_u \to \{ \bp : p(t) \le \bar p(t),  \forall t \in T\}$, with   $\S_u$ as in~\eqref{Eqn_u_alone_policy_space}, using the relation in    \eqref{Eqn_price_prod_relation_cntrl_prblm}: 
$$
 \p(\u) = (\p(u) (t), t \in T ) \quad \mbox{ with } \quad \p(\u)(t) := a u (t) + \frac{\alpha(t)}{2 \beta(t)},
$$
for each $t\in T$. Then 
$J(0, x; \u) := J(0, x; \u, \p(\u) )$  matches with \eqref{Eqn_only_u_dynamics_objective} and $\dot x (t) = u(t) -  ( \alpha(t) - \beta (t) \p(u) (t) )^+ $ matches with the ODE in \eqref{Eqn_main_reduced}.
Thus the proof of  part $(ii)$. \eop

\Forarxiv{\section*{Proof of Theorem \ref{thm_w_dec_inc}}
\label{sec_app_ch3_ch_opt_policy}
  Consider the case when $ x^*(t) > 0 $, for some $t>0$. By hypothesis, 
$$
\lim_{\delta \to 0} \left( w(t -  \delta, t) - w(t , t +  \delta) \right)\le  0. 
$$
If possible assume the contrary---then  further using  the continuity of  $w$ we have a  $\bar \delta >0$, such that 
\begin{equation}
   w(t - \bar \delta, t) > w(t , t + \bar \delta) \mbox{ for all } \delta \le \bar \delta. 
   \label{Eqn_assumption}
\end{equation}
We consider the following modified policy, where     $\epsilon, \delta >0$ will be chosen appropriately (as explained below): set  $  u_{\epsilon}(s)=  u^*(s)$ for $s \notin [t-\delta, t+\delta]$ and otherwise,
\begin{eqnarray*}
    &&u_{\epsilon}(s)  =   
\left ( u^*(s) - \frac{\epsilon}{\gamma(s)}  \right )\mathds{1}_{\{ s \in [t-\delta, t] \} } \\
&& \hspace{30mm}+
  \left ( u^*(s) + \frac{\epsilon}{\gamma(s)} \right )\mathds{1}_{\{  s \in (t, t + \delta] \} }. 
\end{eqnarray*}
By \textbf{B.2},  $u_\epsilon \in {\cal S}_u $, for $\epsilon$   sufficiently small.
Prior to choosing $\epsilon$ and  $\delta$, first     observe   the solution of the ODE from \eqref{Eqn_main_reduced} under $u_\epsilon$ satisfies $x_{\epsilon}(s)= x^*(s)$ if $s \notin [t-\delta, t+\delta]$ and otherwise,
\begin{eqnarray}
&&\hspace{-15mm} x_{\epsilon}(s) = x^*(s) - \epsilon (s- (t-\delta) ) \mathds{1}_{\{ s \in [t - \delta, t] \}} \nonumber\\
&& \hspace{25mm}-  \epsilon \left(\delta-  (s-t) \right)\mathds{1}_{\{s \in [t , t + \delta] \}}.
\label{Eqn_update_policy_traj}
\end{eqnarray}
To begin with, choose $\epsilon, \delta > 0$ such that 
$x_\epsilon(s)  >  0 $ for all $s \in [t-\delta, t+\delta]$: first choose $\delta > 0$  (which  is less than $\bar \delta$) such that ${\underline x} := \inf_{s\in [t-\delta, t+\delta] }x^*(s)  >0$,   possible due to continuity of $x^*$, and then choose $\epsilon >0$ small enough such that  $\epsilon \delta < \underline x$.

From \eqref{Eqn_update_policy_traj}, $x_\epsilon (s) \le x^{*}(s)$ and thus $h(x^{*}(s) ) \ge h(x_\epsilon (s) )$ for all $s$ and $y^{*}(t_1) \ge y_{\epsilon}(t_1)$. And hence, the difference in the objective functions is 
$$
\begin{aligned}
 &J (0, x_{\epsilon};\u_\epsilon)- \rho y_{\epsilon}(t_1) - \left (J(0, x^{*};\u^*) -\rho y^*(t_1)\right ) \\
   &= \int_0^{t_1} \bigg(\frac{\alpha^2(s)}{4 \beta(s)} -a{{{u}_{\epsilon}^*}^2}(s) \gamma(s)-h(x^{*}_{\epsilon}(s))\bigg) ds \\
   & \qquad - \int_0^{t_1} \bigg( \frac{\alpha^2(s)}{4 \beta(s)} -a{u^*}^2(s) \gamma(s)   -h(x^{*}(s))\bigg) ds \\ 
    & \qquad - h_{t_1}(x^{*}_{\epsilon}(t_1)) + h_{t_1}(x^{*}(t_1)) - \rho y_{\epsilon}^*(t_1)  + \rho y^*(t_1) \\
    &=  \int_{t - \delta }^{t} \hspace{-1mm}  -a \bigg(u^{*}(s) - \nicefrac{\epsilon}{\gamma(s)}\bigg)^2  ds +  \int_{t}^{t + \delta } \hspace{-1mm} -a \bigg(u^{*}(s) + \nicefrac{\epsilon}{\gamma(s)}\bigg)^2 ds \\
    &  \qquad  - \int_{t-\delta }^{t + \delta } \bigg( h(x^{*}_{\epsilon}(s)) - h(x^{*}(s))  \bigg) ds + \int_{t-\delta }^{t + \delta }  au^*(s)^2  ds  \\
    & \qquad - h_{t_1}(x^{*}_{\epsilon}(t_1)) + h_{t_1}(x^{*}(t_1))- \rho (y_{\epsilon}^*(t_1)  - y^*(t_1)) \\ 
      &= {2 a \epsilon }\bigg(\int_{t - \delta}^{t}u^*(s) ds - \int_{t}^{t+ \delta}u^*(s) \, ds \bigg) \\
      & \qquad + \int_{t-\delta }^{t + \delta } \bigg( h(x_{\epsilon}(s)) - h(x^{*}(s))  \bigg) \, ds \\
      & \qquad + \int_{t-\delta }^{t + \delta }O(\epsilon^2)  \, ds   - \rho (y_{\epsilon}(t_1)  - y^*(t_1)) >0
\end{aligned}
$$
for some sufficiently small $\epsilon$ in view of 
 \eqref{Eqn_assumption}.  Thus $\u_\epsilon$ performs strictly better contradicting  the fact that the policy $\u^{*}$ is optimal. 
 
 With almost similar arguments, one can prove part $(ii)$ corresponding to  $ x^{*}(t) < 0$. The main differences to be observed  here are: a) we will have $0 > x_\epsilon(s) \ge  x^{*}(s) $ in the selected neighbourhood of $t$; b) however the magnitude is smaller, i.e.,  $x^2_\epsilon(s) \le {x^{*}}^2(s) $, and so $h(x^{*}(s) ) \ge h(x_\epsilon(s))$ again for all $s$; and    c)  since $x^{*}(0) = x_\epsilon(0) \ge  0$, we have $y(t_1) = y_{\epsilon}(t_1)$, where $y(t_1) := \sup_{s \in [0,t_1]} \{x(s)\}$ and $ y_{\epsilon}(t_1) := \sup_{s \in [0,t_1]} \{x_{\epsilon}(s)\}$. 
 
\textit{Proof of part (ii):}   If possible, consider $x^{*}(t_1) > 0$. In this case, we introduce a modified policy by setting $\epsilon, \delta >0$, such that the policy $\u_{\epsilon}$ and the corresponding  $\x_{\epsilon}$ is
\begin{eqnarray*}
u_{\epsilon}(s)  &=&  \left \{  \begin{array}{lll}
 u^*(s) - \epsilon       &  \mbox{if }   s \in [t_1 -\delta, t_1]  \\
   u^*(s)    & \mbox{else }, 
\end{array} \right. \mbox{ and, } \\
x_{\epsilon}(s)  &=& \left \{  \begin{array}{lll}
 x^*(s) - \epsilon (s- (t_1 -\delta) )        &  \mbox{if }   s \in [t_1 - \delta, t_1 ]  \\
   x^*(s)    & \mbox{else}.
\end{array} \right.
\end{eqnarray*}
To ensure $x_{\epsilon}(t_1) $ remains positive, we carefully select $\epsilon, \delta $.  Then by using the same arguments as in part (i), difference in the objective function $J (0, x_\epsilon; \u_\epsilon) - J(0, x^*; \u^*)$ is strictly positive.
 By carefully selecting $\delta > 0 $, we can demonstrate an improvement in the objective function. Thus, we conclude that $ x^*(t_1) \leq 0 $. \eop
}

\section*{Proof of Proposition \ref{thm_y_terminal_dec_with_sigma}} 
\label{sec_app_ch3_thm_y_terminal_dec_with_sigma}

Consider $\tilde{\rho} < \rho$. Let $\x_{\tilde{\rho}}^*, \u_{\tilde{\rho}}^*$ and $\x_\rho^*, \u_\rho^*$ be the optimal trajectory and optimal policy wrt $\Tilde{\rho}$ and $\rho$ respectively. Define 
 $$
       \begin{aligned}
            J^{*}_{o} &:= \int_0^{t_1} \bigg( \frac{\alpha^2(s)}{4 \beta(s)} -a {u_\rho^*}^2(s) \gamma(s)   -h(x_\rho^*(s))\bigg) ds - h_{t_1}(x_\rho^*(t_1)),  \\
        \tilde{J}^{*}_{o} &:= \int_0^{t_1} \bigg( \frac{\alpha^2(s)}{4 \beta(s)} -a {u_{\tilde{\rho}}^*}^2(s) \gamma(s)  -h(x_{\tilde{\rho}}^*(s))\bigg) ds  - h_{t_1}(x_{\tilde{\rho}}^*(t_1)).
       \end{aligned}
       $$
       To begin with,  observe  the value function $v_\rho = J_o^* - {\rho} y_{{\rho}}^*(t_1) \ge  \tilde J_o^* - { \rho} y_{\tilde{\rho}}^*(t_1) $, the objective function corresponding to $\rho$ and at policy  $\tilde \u^*_\rho$.  
        Assume, if possible, $y_{{\rho}}^*(t_1) > y_{\tilde{\rho}}^*(t_1) $, then:
\begin{eqnarray*}
J_o^* - \tilde{\rho} y_{{\rho}}^*(t_1)    &=& J_o^* - {\rho} y_{{\rho}}^*(t_1) + (\rho -\tilde{\rho} )y_{{\rho}}^*(t_1) \\ 
 &\geq&   \tilde{J}^{*}_{o} - {\rho} y_{\tilde{\rho}}^*(t_1) + (\rho -\tilde{\rho} ) y_{{\rho}}^*(t_1) \\ 
     & =&   \tilde{J}^{*}_{o} - \tilde{\rho} y_{\tilde{\rho}}^*(t_1) + (\rho -\tilde{\rho} ) \bigg(y_{{\rho}}^*(t_1) - y_{\tilde{\rho}}^*(t_1)\bigg)  \\
    &>&   \tilde{J}^{*}_{o} - \tilde{\rho} y_{\tilde{\rho}}^*(t_1) = v_{\tilde \rho} .
\end{eqnarray*}
The last inequality arises from $ \tilde{\rho} < \rho $ and $ y_{{\rho}}^*(t_1) > y_{\tilde{\rho}}^*(t_1) $. This leads to a contradiction to the optimality  of $  \u^*_{\tilde \rho}$ for the problem with parameter $\tilde \rho$. \eop

\end{document}


\maketitle

\section{A detailed example}

Here we include some equations and theorem-like environments to show
how these are labeled in a supplement and can be referenced from the
main text.
Consider the following equation:
\begin{equation}
  \label{eq:suppa}
  a^2 + b^2 = c^2.
\end{equation}
You can also reference equations such as \cref{eq:matrices,eq:bb} 
from the main article in this supplement.

\lipsum[100-101]

\begin{theorem}
An example theorem.
\end{theorem}

\lipsum[102]
 
\begin{lemma}
An example lemma.
\end{lemma}

\lipsum[103-105]

Here is an example citation: \cite{KoMa14}.

\section[Proof of Thm]{Proof of \cref{thm:bigthm}}
\label{sec:proof}

\lipsum[106-112]

\section{Additional experimental results}
\Cref{tab:foo} shows additional
supporting evidence. 

\begin{table}[htbp]
\footnotesize
  \caption{Example table.}  \label{tab:smfoo}
\begin{center}
  \begin{tabular}{|c|c|c|} \hline
   Species & \bf Mean & \bf Std.~Dev. \\ \hline
    1 & 3.4 & 1.2 \\
    2 & 5.4 & 0.6 \\ \hline
  \end{tabular}
\end{center}
\end{table}

\bibliographystyle{siamplain}
\bibliography{references}